\documentclass[a4paper,11pt]{amsart}
\pdfoutput=1

\usepackage{amssymb}
\usepackage{graphics}
\usepackage[english]{babel}
\usepackage{ifthen}
\usepackage{ifpdf}
\usepackage{enumerate,amsfonts,amscd,amsthm,amssymb,cite} 
\usepackage{graphicx}
\usepackage{url}
\usepackage{todonotes}

    \usepackage[  plainpages = false, pdfstartview = {FitH},
		 pdfpagelabels, pdfpagelayout = SinglePage,
		 bookmarks = true, bookmarksopen = true,
     bookmarksnumbered = true, linktocpage,
     colorlinks = true,
     linkcolor = blue, urlcolor  = blue,
     citecolor = red, anchorcolor = green,
     hyperindex = true, hyperfigures]{hyperref} 

\def\ds{\displaystyle}

\newcommand{\maxo}{\overline{\omega}}
\newcommand{\mino}{\underline{\omega}}

\newcommand{\gotable}[3]
      {\begin{table}[h] 
       \begin{center}
       \addtolength{\tabcolsep}{#1mm}
       \renewcommand{\arraystretch}{#2}
       \begin{tabular}{#3}}

\newcommand{\stoptable}
      {\end{tabular}
       \end{center}
       \end{table}}

\newcommand{\beq} {\begin{eqnarray}}
\newcommand{\eeq}{\end{eqnarray}}

\numberwithin{equation}{section}

\theoremstyle{plain}
\newtheorem{theorem}{Theorem}[section]
\newtheorem{lem}[theorem]{Lemma}
\newtheorem{cor}[theorem]{Corollary}
\newtheorem{prop}[theorem]{Proposition}

\theoremstyle{definition}
\newtheorem{defn}[theorem]{Definition}
\newtheorem{ex}[theorem]{Example}
\newtheorem{xca}[theorem]{Exercise}
\newtheorem{conj}[theorem]{Conjecture}

\theoremstyle{remark}
\newtheorem{rem}[theorem]{Remark}

\newcommand{\bt}[1]{\begin{theorem}\label{#1}}
\newcommand{\bc}[1]{\begin{cor}\label{#1}}
\newcommand{\bl}[1]{\begin{lem}\label{#1}}
\newcommand{\bp}[1]{\begin{prop}\label{#1}}
\newcommand{\be}[1]{\begin{ex}\label{#1}}
\newcommand{\bd}[1]{\begin{defn}\label{#1}}
\newcommand{\br}[1]{\begin{rem}\label{#1}}
\newcommand{\bx}[1]{\begin{xca}\label{#1}}
\newcommand{\bcon}[1]{\begin{conj}\label{#1}}

\newcommand{\et}{\end{theorem}}
\newcommand{\ec}{\end{cor}}
\newcommand{\el}{\end{lem}}
\newcommand{\ep}{\end{prop}}
\newcommand{\ee}{\end{ex}}
\newcommand{\ed}{\end{defn}}
\newcommand{\exc}{\end{xca}}
\newcommand{\er}{\end{rem}}
\newcommand{\econ}{\end{conj}}

\newcommand{\bpr}{\begin{proof}}
\newcommand{\epr}{\end{proof}}
\def\A  {\mathcal{A}}
\def\FA {\mathcal{F(A)}}
\def\Ch {\mathcal{C(\A)}}

\def \F {\mathcal{F}}

\def \P {\mathcal{P}}

\def \G {\mathcal{G}}

\def \L {\mathcal{L}}
\def \T {\mathcal{T}}
\def \h{\mathcal{H}}

\def\R{\mathbb{R}}

\def\C{\mathbb{C}}

\title{On Arrangements of Pseudohyperplanes}
\author{Priyavrat Deshpande}
\address{Chennai Mathematical Institute\\ H1 SIPCOT IT park, Sirueri\\ Tamil Nadu, India}
\email{pdeshpande@cmi.ac.in}
\date{}

\begin{document}

\begin{abstract} To every realizable oriented matroid there corresponds an arrangement of real hyperplanes. 
The homeomorphism type of the complexified complement of such an arrangement is completely determined by the oriented matroid. 
In this paper we study arrangements of pseudohyperplanes; they correspond to non-realizable oriented matroids. 
These arrangements arise as a consequence of the Folkman-Lawrence topological representation theorem. 
We propose a generalization of the complexification process in this context. 
In particular we construct a space naturally associated with these pseudo-arrangements which is homeomorphic to the complexified complement in the realizable case. 
Further we generalize the classical theorem of Salvetti and show that this space has the homotopy type of a cell complex defined in terms of the oriented matroid.\end{abstract}
\maketitle

\section*{Introduction}

An arrangement of hyperplanes is a finite set $\A$ consisting of affine hyperplanes in $\R^l$. These hyperplanes and their intersections induce a stratification of $\R^l$. 
These strata, called the \textit{faces of $\A$}, form a poset ordered by topological inclusion and the set of all possible intersections forms a poset ordered by reverse inclusion.
These posets contain combinatorial information about the arrangement. An important topological object associated with an arrangement $\A$ is the \emph{complexified complement} $M(\A)$. 
It is the complement of the union of the complexified hyperplanes in $\C^l$. 
One of the important aspects of the theory of arrangements is to understand the interaction between the combinatorial data of an arrangement and the topology of this complement. 
For example, a pioneering result by Salvetti \cite{sal1} states that the homotopy type of the complement is determined by the face poset. 
He constructed a regular cell complex (now known as the \textit{Salvetti complex}), using the incidence relations among faces, on which the complement deformation retracts (Bj\"orner and Ziegler \cite{bz92} extended this result to show that the homeomorphism type of the complement is also determined by the face poset).

\par

Oriented matroids are intimately connected with central hyperplane arrangements (i.e., all hyperplanes are passing through the origin). 
This combinatorial structure not only combines the above mentioned posets but also supplies rich techniques to study various aspects of arrangements. 
The strata of a central hyperplane arrangement satisfy the covector axioms of an oriented matroid. 
The oriented matroids which correspond to faces of a hyperplane arrangement are known as the \emph{realizable oriented matroids}. There are oriented matroids that do not correspond to hyperplane arrangements (e.g., non-Pappus configuration, see Example \ref{ex1s1c5}). 
Hence for a long time an important question in this field was to come up with the right topological model for oriented matroids. 
This was settled by Folkman and Lawrence in \cite{fl78}. The \emph{Folkman-Lawrence Topological Representation Theorem} states that in general oriented matroids correspond to certain collections of finitely many topological spheres
and balls (i.e. \emph{arrangements of pseudo-hemispheres}). 
These so-called pseudo-arrangements not only describe oriented matroids in the same way that $\R^l$ and collections of half spaces describe an obvious combinatorial structure but there is a one-to-one correspondence between such arrangements and the oriented matroids. 
In his thesis Mandel \cite{mandel82} introduced \emph{sphere systems} that simplified some aspects of the pseudo-hemisphere arrangements and also proved the piecewise-linear version of the representation theorem. \par

In his thesis Ziegler \cite[Section 5.5]{zieg87} extended the definition of the Salvetti complex to arbitrary oriented matroids. To every oriented matroid one can associate a simplicial complex and in case of a realizable oriented matroid this complex has the homotopy type of the space $M(\A)$. In their paper Gel'fand and Rybnikov \cite{gr89} studied the Salvetti complex for arbitrary oriented matroids and showed that the cohomology ring of this complex is isomorphic to the Orlik-Solomon algebra of the associated lattice of flats (see also \cite{bz92}). This result not only extends the classical theorem of Brieskorn and Orlik-Solomon but also gives a completely combinatorial proof. \par

An important object missing in the context of pseudo-arrangements is an analogue of the space $M(\A)$. 
The aim of this paper is to introduce such a space.
We represent oriented matroids using arrangements of pseudo-hyperplanes; these are collections of subspaces obtained by mildly deforming linear hyperplanes.
The intersection of such an arrangement with the unit sphere gives the pseudo-sphere arrangement associated with the oriented matroid.
Since these pseudo-hyperplanes are not necessarily defined by vanishing of polynomials the usual complexification does not work.
We make use of the theory of microbundles introduced by Milnor \cite{milnor64} to get around this difficulty.
Note that the total space of the tangent microbundle of $\R^l$ is homeomorphic to $\C^l$ and that of a pseudo-hyperplane is homeomorphic to $\C^{l-1}$.
Hence as an analogue of the space $M(\A)$ we consider the complement of the union of tangent microbundles of the pseudo-hyperplanes inside $\R^{2l}$.
In case the arrangement is a collection of honest hyperplanes then our construction yields a space homeomorphic to the complexified complement.
\par 

The paper is organized as follows. We begin in Section \ref{sec1} by recalling relevant results about hyperplane arrangements. We also recall the basic notions from oriented matroids and link them to hyperplane arrangements. In Section \ref{sec2} we first state the topological representation theorem. For the convenience of the reader we give a slightly detailed exposition. Then we introduce the generalization of the complexified complement and prove our main result. Section \ref{section3} is about generalizing the relationship between hyperplane arrangements and zonotopes. We explore metrical-hemisphere complexes as a possible generalization. Finally in Section \ref{section4} we gather results from the theory of hyperplane arrangements which generalize in this context without much work. We also discuss the relevance of our work and outline future research. 

\subsection*{Acknowledgments} This paper is a part of the author's doctoral thesis \cite[Chapter 5]{deshpande_thesis11}. 
The author would like to thank his supervisor Graham Denham for his support. 
The author would also like to thank Eric Babson, Emanuele Delucchi, Alex Papazoglou and Thomas Zaslavsky for fruitful discussions. 
The author would like to acknowledge the support of the Mathematics department at Northeastern University for hosting a visit during 2011-12 during which the initial draft was written. 
Sincere thanks to the anonymous referee for many helpful suggestions and for pointing out the connection with microbundles.

\section{Arrangements and Oriented Matroids}\label{sec1}
We start by briefly reviewing some basic facts. 

\subsection{Hyperplane Arrangements} \label{sec:prelim}

Hyperplane arrangements arise naturally in geometric, algebraic and combinatorial instances. 
In this section we define hyperplane arrangements and the associated combinatorial and topological data in a setting that is most relevant to our work. 

\bd{def1} A real \emph{arrangement of hyperplanes} is a collection $\A = \{H_1,\dots,H_k\}$ of finitely many \emph{hyperplanes} in $\R^l$, $l\geq 1$. \ed

The \textit{rank of an arrangement} is the largest dimension of the subspace spanned by the normals to the hyperplanes in $\A$. We assume that all our hyperplane arrangements are \textit{essential}, i.e., their rank is equal to the dimension of the ambient vector space. An arrangement is called \textit{central} if the intersection of all the hyperplanes in $\A$ is non-empty. For a subset $X$ of $\R^l$, the \emph{localization} of $\A$ at $X$ is  the sub-arrangement $\A_X := \{H\in\A ~|~ X\subseteq H\}$.  \par

There are two posets associated with $\A$, namely, the intersection poset and the face poset which contain the combinatorial information. 

\bd{def2} The \emph{intersection poset} $L(\A)$ of $\A$ is defined as the set of all intersections of hyperplanes ordered by reverse inclusion. \ed

$L(\A)$ is a ranked poset with the rank of an element being the codimension of the corresponding intersection. In general it is a (meet) semilattice; it is a lattice if and only if the arrangement is central. \par 

The hyperplanes of $\A$ induce a stratification of $\R^l$ such that each stratum is an open polyhedron; these strata are called the \emph{faces} of $\A$.

\bd{def3} The \emph{face poset} $\FA$ of $\A$ is the set of all faces ordered by topological inclusion: $F\leq G$ if and only if $F\subseteq\overline{G}$. \ed

The codimension-$0$ faces are called \emph{chambers}. The set of all chambers is denoted by $\Ch$. A chamber is \textit{bounded} if it is a bounded subset of $\R^l$. Two chambers $C$ and $D$ are \textit{adjacent} if their closures have a non-empty intersection containing a codimension-$1$ face.
By a \textit{complexified hyperplane arrangement} we mean an arrangement of hyperplanes in $\C^l$ for which the defining equation of each hyperplane is real. Hence to every hyperplane arrangement in $\R^l$ there corresponds an arrangement of hyperplanes in $\C^l$. An important topological space associated with a real hyperplane arrangement is the following.

\bd{def4} Let $\A$ denote a hyperplane arrangement in $\R^l$. The \emph{complexified complement} of $\A$ is denoted by $M(\A)$ and is defined as 
\[M(\A) := \C^l \setminus \ds (\bigcup_{H\in\A} H_{\C}) \] 
where $H_{\C}$ is the hyperplane in $\C^l$ with the same defining equation as $H$.\ed 

Note that since $M(\A)$ is the complement of subspace of real codimension $2$ in $\C^l$, it is connected. It is an open submanifold of $\C^l$ with the homotopy type of a finite-dimensional CW complex \cite[Section 5.1]{orlik92}.


We explain the construction of the regular $l$-complex, called the \textit{Salvetti complex} and denoted by $Sal(\A)$, by first describing its cells. The $k$-cells, for $0\leq k\leq l$, of $Sal(\A)$ are in one-to-one correspondence with the pairs $[F, C]$, where $F$ is a codimension-$k$ face of $\A$ and $C$ is a chamber whose closure contains $F$.
 
We specify the boundary of each cell - by Exercise 2.21 and Proposition 4.7.23 of \cite{ombook99} these data uniquely determine a regular cell structure. 
A cell labeled $[F_1, C_1]$ is contained in the boundary of another cell labeled $[F_2, C_2]$ if and only if $F_2 \leq F_1$ in $\FA$ and $C_1, C_2$ are contained in the same chamber of $\A_{F_1}$. \par 

The Salvetti complex has the same homotopy type as that of the complexified complement.

\begin{theorem}[Salvetti \cite{sal1}]\label{thm0} Let $\A$ be an arrangement of real hyperplanes and $M(\A)$ be the complement of its complexification inside $\C^l$. Then there is an embedding of $Sal(\A)$ into $M(\A)$. Moreover there is a natural map in the other direction which is a deformation retraction.\et

The cohomology algebra of $M(\A)$, called the \emph{Orlik-Solomon algebra}, is particularly interesting. Its construction is completely combinatorial and depends on the data that are encoded by the intersection poset of $\A$ (see \cite[Chapter 3, Section 5.4]{orlik92} for details). 

\subsection{Oriented Matroids}\label{ch1secOM} 
Let $E = \{1,\dots, n\}$ be the finite ground set for some $n>0$. A \emph{sign vector} is a function $X\colon E \to \{+, 0, -\}$, i.e., an assignment of signs to the elements of $E$. 
The set of all possible sign vectors is denoted by $\{ +, 0, -\}^E$ and $X_e$ stands for $X(e)$ for all $e\in E$. 
The \emph{support} of a vector $X$ is $\underline{X} = \{e\in E | X_e\neq 0\}$; its \emph{zero set} is $z(X) = E\setminus \underline{X}$. 
The opposite of a vector $X$ is $-X$, defined by $(-X)_e = -(X_e)$. The \emph{zero} vector is $\textbf{0} = (0,\dots,0)$.
The \emph{composition} of two sign vectors $X$ and $Y$ is $X\circ Y$ defined by 
\[(X\circ Y)_e := \begin{cases} X_e, &\hbox{if~} X_e\neq 0, \\
Y_e, &\hbox{otherwise}.\end{cases}\]
The \emph{separation set} of $X$ and $Y$ is $S(X, Y) = \{e\in E | X_e = -Y_e \neq 0\}$. With these terminologies we can define oriented matroids using the \emph{covector axioms}. 

\bd{def16} A set $\L\subseteq \{-, 0, +\}^E$ of signed vectors is the set of covectors of an oriented matroid if it satisfies: 
\begin{enumerate}
	\item[(V0)] $\textbf{0}\in \L$,
	\item[(V1)] $X\in \L\Rightarrow -X\in \L$,
	\item[(V2)] $X, Y\in \L\Rightarrow X\circ Y \in \L$,
	\item[(V3)] if $X, Y\in \L$ and $e\in S(X, Y)$ then there exists $Z\in \L$ such that $Z_e = 0$ and $Z_f = (X\circ Y)_f = (Y\circ X)_f$ for all $f\notin S(X, Y)$. 
\end{enumerate} \ed

There is a partial order on the sign vectors defined as follows:
\[Y\leq X \iff Y_e \in \{0, X_e \} \quad\forall~ e\in E.\]

If $\L \subset \{+, 0, -\}^E$ is a set of covectors of an oriented matroid then it inherits the above defined partial ordering to become a poset with the bottom element $\textbf{0}$. 
The poset $\hat{\L} := (\L\cup \{\hat{1}\}, \leq)$ is a lattice. 
The join  of $X$ and $Y$ is $X\circ Y = Y\circ X$ if $S(X, Y) = \emptyset$, and equals $\hat{1}$ otherwise. 

\bd{def2s1c5}
The lattice $\F(\L) = (\hat{\L}, \leq)$ is called the \emph{face lattice} of the oriented matroid $\L$. 
The maximal elements of $\L$ are called \emph{topes} (or \emph{regions}).
Let $\T(\L)$ denote the set of topes. 
The \emph{rank} of $\L$ is the length of a maximal chain in $(\L, \leq)$.\ed

We describe how oriented matroids arise in the context of central hyperplane arrangements. 
Let $\A =\{H_1,\dots,H_n \}$ be a central arrangement of hyperplanes in $\R^l$. 
Associated to every hyperplane $H_i\in \A$, there are two open half-spaces bounded by the hyperplane, which will be denoted by $H_i^+$ and $H_i^-$. 
Note that the linear form defining the hyperplane $H_i$ maps the points in $H^+_i$ side onto $\R_{>0}$ and the points in $H_i^-$ onto $\R_{<0}$. 
We assign a sign vector $X(v) = (X_1(v),\dots,X_n(v))$ to every point $v\in \R^l$ as follows: 
\[ 
	X_i(v) = \begin{cases} + &\hbox{if~} x\in H_i^+, \\
					0 &\hbox{if~} x\in H_i, \\
					- &\hbox{if~} x\in H_i^-. \end{cases} \]
					
Let $\F$ denote the set of all possible sign vectors that arise due to the induced stratification. 
It is clear that the faces of $\A$ are in one-to-one correspondence with the sign vectors. 
They satisfy the above mentioned axioms for oriented matroids.
The face poset of $\A$ is isomorphic to the poset of covectors of the oriented matroid.
Hence every central hyperplane arrangement gives rise to an oriented matroid.
An oriented matroid which arises in this way is called a \emph{realizable oriented matroid}.\par 



\br{rem1s1c5} If $\L$ is a (linear) oriented matroid arising from a central hyperplane arrangement $\A$ in $\R^l$ then $\F(\L)$ is isomorphic to the face poset $\F(\A)$. 
In particular the topes of $\L$ correspond to the chambers of $\A$. 
The composition of two covectors and their separation set also have a geometric meaning which we now state.
Say that a hyperplane separates two chambers if these chambers lie in the distinct connect components of the complement of this hyperplane.
If $X, Y$ and are two topes then $S(X, Y)$ corresponds bijectively to the set of hyperplanes separating the corresponding chambers.
For any two chambers define the distance between them as the number of separating hyperlanes.
Note that given a face $F$ and chamber $C$ there exists a unique chamber $C_F$ which is incident with $F$ and closest to $C$ with respect to the above distance.
If $X$ corresponds to a $F$ and $Y$ corresponds to $C$ then the composition $X\circ Y$ corresponds to the chamber $C_F$.
\er  



\section{pseudo-hyperplane Arrangements}\label{sec2}
An arbitrary oriented matroid need not correspond to an arrangement of hyperplanes. Such oriented matroids are called \emph{non-realizable oriented matroids}. The \emph{Folkman-Lawrence topological representation theorem} \cite{fl78} states that every oriented matroid is `\emph{almost}' realizable. Originally the topological representation theorem was stated in terms of \emph{pseudo-hemisphere arrangements}. Later Arnaldo Mandel in his thesis \cite{mandel82} achieved much simplification. He reproved the theorem in terms of \emph{sphere systems} (also known as \emph{arrangements of pseudospheres}). \par

This section begins with a brief review of the topological representation theorem. We assume reader's familiarity with the basic PL topology; a relevant reference, in this context, is Rushing's book \cite[Chapter 1]{rushing73}. For basic concepts in topological combinatorics like geometric realization of a poset, simplicial complexes, nerve lemma etc. we refer the reader to Kozlov's book \cite{dk1}. Our review of the topological representation theorem is based on \cite[Chapter 5]{ombook99}. We reproduce it here for the benefit of the reader. After this review we introduce a generalization of hyperplane arrangements called \textit{arrangements of pseudo-hyperplanes}.\par

The pseudo-hyperplanes are obtained by mildly deforming hyperplanes. 
We will show that their arrangements correspond to pseudosphere arrangements and hence to oriented matroids. 
Since the pseudo-hyperplanes may not be described by algebraic equations in general, it is not possible to define their complexification. 
For each pseudo-hyperplane arrangement in $\R^l$ we construct a subspace of $\R^{2l}$ which plays the role of the complexified complement. 
The main result of this section is the proof that this subspace has the homotopy type of the associated Salvetti complex. 

\subsection{The topological representation theorem}\label{sec2p1} 
An $n$-manifold $N$ contained in an $l$-manifold $X$ is \emph{locally flat} at $x\in N$ if there exists a neighborhood $U_x$ of $x$ in $X$ such that $(U_x, U_x\cap N)\cong (\R^l, \R^n)$. An embedding $f\colon\thinspace N\to X$ such that $f(N)\subseteq X$ is said to be \emph{locally flat at a point} $x\in N$ if $f(N)$ is locally flat at $f(x)$. Embeddings and submanifolds are \emph{locally flat} if they are locally flat at every point. Here it is assumed that the manifolds and submanifolds we consider are without boundary.

\begin{theorem}[{\cite[Theorem 1.7.2]{rushing73}}]\label{t2ref} Let $S^l$ denote the standard unit sphere in $\R^{l+1}$. If $f\colon S^l\to \R^n$, $n-l\neq 2$ is a locally flat embedding then there exists a homeomorphism $h\colon \R^n\to \R^n$ such that $h\circ f$ is the inclusion map. The same conclusion holds for an embedding of $\R^l$ into $\R^n$. \et

A subset of the standard unit sphere is called a subsphere if it is homeomorphic to some lower-dimensional sphere. We single out a class of subspheres that play an important role in defining more general types of arrangements. The following result is from \cite[Theorem 1.8.2]{rushing73}.

\bl{lem1s1c5} For an $(l-1)$-subsphere $S$ of $S^l$ the following conditions are equivalent:
\begin{enumerate}
	\item the embedding of $S$ is equivalent to the inclusion of the standard $(l-1)$-sphere in $S^l$;
	\item the embedding of $S$ is equivalent to some PL $(l-1)$-subsphere of $S^l$;
	\item the closure of each connected component of $S^l\setminus S$ is homeomorphic to the $l$-ball. \end{enumerate} \el 

The equivalence class of these (embeddings of) subspheres is known as \emph{tame}, all other embeddings are called \emph{wild}. It is known that all embeddings of $S^1$ into $S^2$ are tame (the Sch\"onflies theorem \cite{rushing73}). However, there are wild $2$-spheres in $S^3$, for example, the Alexander horned sphere (see \cite[Page 69]{rushing73} and \cite[Example 2B.2]{hatcher02}). 
	
\bd{def3s1c5}
An $(l-1)$-subsphere $S$ in $S^l$ satisfying any of the equivalent conditions in Lemma \ref{lem1s1c5} is called a \emph{pseudosphere} in $S^l$. The two connected components of $S^l\setminus S$ are its \emph{sides}, denoted by $S^+$ and $S^-$. The closures of the sides are called the \emph{closed sides} (or \emph{pseudohemispheres}).\ed

We can now introduce the pseudo-arrangements that were used to prove the representation theorem. 

\bd{def4s1c5}
A \emph{signed arrangement of pseudospheres} in the standard unit sphere $S^l\subseteq \R^{l+1}$ is a finite collection $\A = \{(S^+_i, S^0_i, S^-_i)~|~ i\in E \}$ where $E = \{1,\dots, n \}$ such that 
\begin{enumerate}
	\item each $S^0_i$ is a pseudosphere in $S^l$ with sides $S^+_i$ and $S^-_i$,
	\item if for some subset $I$ and an index $j$ the intersection $S_I := \bigcap_{i\in I} S_i^0\nsubseteq S_j^0$, then $S_I\cap S_j^0$ is a pseudosphere in $S_I$ with sides $S_I\cap S^+_j$ and $S_I\cap S^-_j$. \end{enumerate}\ed

Such an arrangement is said to be \textit{essential} if the intersection of all the pseudospheres is empty. Unless otherwise stated all arrangements are signed and essential. Moreover, for the sake of notational simplicity we assume that both the sides of each pseudosphere are equipped with a sign and we will not explicitly mention it every time. Since each side has a sign attached to it one can define a sign function similar to that for hyperplane arrangements. The position of each point $x\in S^l$ with respect to each pseudosphere in the arrangement $\A$ is given by a sign vector $\sigma(x)\in \{+, 0, -\}^E$, defined by
\[ 
	\sigma(x)_i = \begin{cases} + &\hbox{if~} x\in S_i^+, \\
					0 &\hbox{if~} x\in S_i^0, \\
					- &\hbox{if~} x\in S_i^-. \end{cases} \]
					
The reader can verify that such an arrangement defines a stratification of the ambient sphere. We refer to these strata as \textit{faces} of that arrangement. Each face is indexed by a sign vector in $\sigma(S^l)$. One of the reasons for this generalization is the following:

\bt{thm1s1c5} Let $\A$ be an arrangement of pseudospheres in $S^l$. Then 
\[\L(\A) := \{\sigma(x)~|~ x\in S^l \}\cup \{\textbf{0}\}\subseteq \{+, 0, -\}^E \]
is the set of covectors of an oriented matroid and the rank of $\L(\A)$ is $l+1$. \et


\begin{lem}[{\cite[Lemma 3, page 201]{mandel82}}]\label{lem2s1c5} Let $\A$ be a signed and essential arrangement of pseudospheres in $S^l$. For every $X\in \L(\A)\setminus \{\textbf{0}\}$ the face $\sigma^{-1}(X)$ is an open cell of a regular cell decomposition $\Delta(\A)$ of $S^l$. The boundary of $\overline{\sigma^{-1}(X)}$ is the union of all those $\sigma^{-1}(Y)$ such that $Y < X$. Furthermore, the mapping $X\mapsto \{y\in S^l~|~ \sigma(y)\leq X\}$ gives an isomorphism 
\[\hat{\L}(\A) \cong \hat{\F}(\Delta(\A)) \] of the face lattice of $\L(\A)$ and face lattice of the regular cell complex $\Delta(\A)$. \el

Two signed arrangements $\A = \{S_1, \dots, S_n\}$ and $\A' = \{S'_1,\dots, S'_n \}$ of pseudospheres in $S^l$ are \emph{topologically equivalent}, denoted $\A \sim \A'$, if there exists a homeomorphism $h\colon S^l\to S^l$ such that $h(S_i) = S'_i$ and $h(S^+_i) = (S'_i)^+$ for all $1\leq i\leq n$. This topological equivalence is combinatorially determined. 

\bt{thm2s1c5} Two signed arrangements $\A$ and $\A'$ in $S^l$ are topologically equivalent if and only if $\L(\A) \cong \L(\A')$. \et


\bc{cor1s1c5} Let $\A = \{S_1, \dots, S_n \}$ be a signed arrangement of pseudospheres in $S^l$. The oriented matroid $\L(\A)$ is realizable if and only if there exists a homeomorphism $h\colon S^l\to S^l$ such that $h(S_i) = S^l\cap H_i$, where $H_i$ is a codimension-$1$ linear subspace of $\R^{l+1}$, for every $i$.  \ec

An arrangement is said to be \emph{centrally symmetric} if each pseudosphere $S\in \A$ is invariant under the antipodal mapping of $S^l$ (and so are the sides, i.e., $S_i^+\mapsto S_i^-$ for every $i$). 

\begin{theorem}[{Topological representation theorem \cite[Theorem 5.2.1]{ombook99}}]\label{thm3s1c5} Let $\L\subseteq \{+, 0, -\}^E$. Then following conditions are equivalent:
\begin{enumerate}
	\item $\L$ is the set of covectors of a (simple) oriented matroid of rank $l$. 
	\item $\L = \L(\A)$ for some signed arrangement $\A = \{S_1,\dots, S_n\}$ of pseudospheres in $S^{l-1}$, which is essential and centrally symmetric and whose induced cell complex $\Delta(\A)$ is regular.
\end{enumerate} \et

Let $\L$ be the set of covectors of a rank $l$ oriented matroid. According to Theorem \ref{thm3s1c5} there corresponds a signed arrangement $\A = \{S_1,\dots, S_n \}$ of pseudospheres in $S^{l-1}$, the unit sphere in $\R^l$. Since each pseudosphere $S$ is centrally symmetric any pair of antipodal points $x, -x\in S$ generates a line through the origin in $\R^l$. For $S\in \A$ let $H_S$ be the set of all rays from the origin passing through $S$. In fact $H_S$ is the cone over $S$; we call it a \textit{pseudo-hyperplane}.
The next result is now immediate and follows from Lemma \ref{lem1s1c5}. 

\bl{lem3c1s5}
Let $S$ be a pseudosphere in the unit sphere $S^{l-1}$ and $H_S$ be the cone. Then there exists a homeomorphism of $\R^l$ such that it maps $H_S$ to a hyperplane passing through the origin. \el

\bd{def5c1s5} An \emph{arrangement of pseudo-hyperplanes} is a finite collection $\A$ of pseudo-hyperplanes in $\R^l$ such that $\{H\cap S^{l-1}~|~H\in\A\}$ is an arrangement of pseudospheres in $S^{l-1}$. \ed

Given an arrangement $\A$ of pseudospheres in $S^{l-1}$ we denote by $c\A$ the corresponding arrangement of pseudo-hyperplanes. A nonzero face of $c\A$ is the cone over some face of $\A$ and hence homeomorphic to an open polyhedral cone of one dimension higher. As before we denote the set of all faces of $c\A$ by $\F(c\A)$ and the set of all chambers, i.e., the top-dimensional faces by $\mathcal{C}(c\A)$. We suppress the reference to $c\A$ if the context is clear. 

\be{ex1s1c5}
Consider the non-Pappus arrangement of $9$ pseudolines in $\R^2$. Corresponding to this arrangement there is a rank $3$ oriented matroid on $10$ elements which realizes an arrangement of pseudocircles in $S^2\subseteq \R^3$, see Figure \ref{nonpappus} below. we refer the reader to \cite[Theorem 7.26, Example 7.28]{zig95} for a detailed explanation of why this arrangement is not realizable and a description of the corresponding oriented matroid.
 \begin{figure}[!ht]
  \begin{center}  
    \includegraphics[scale=0.2,clip]{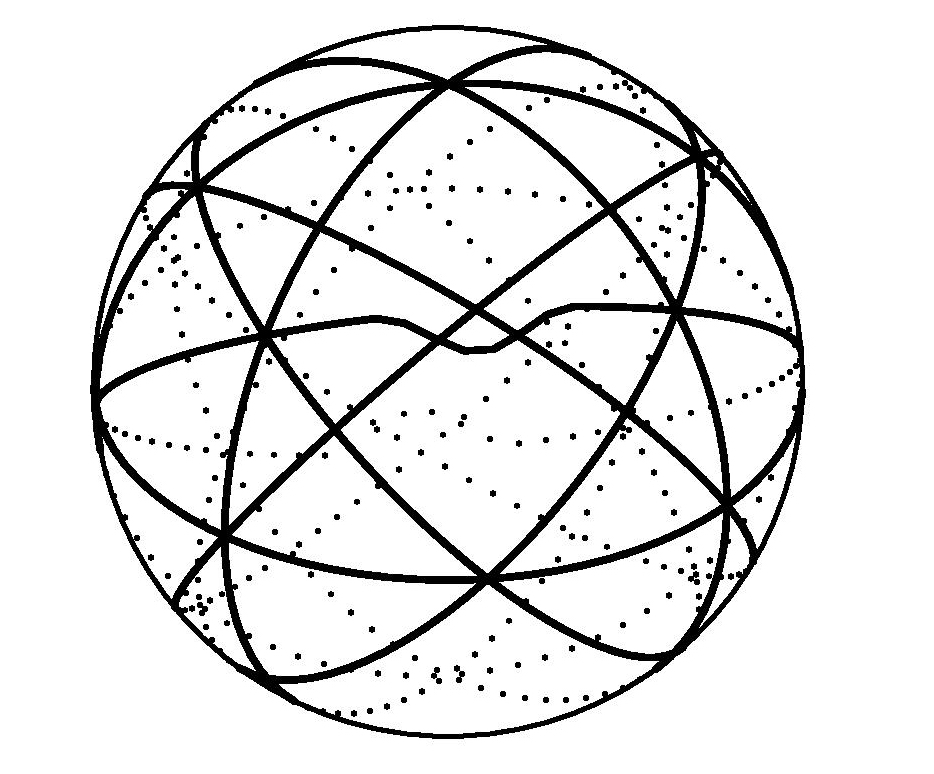} 
  \end{center}
  \caption{Non-Pappus pseudocircle arrangement.} \label{nonpappus}     
 \end{figure}
Now consider cones over these pesudocircles. The cones are pseudoplanes in $\R^3$. Each pseudoplane, individually, is homeomorphic to a $2$-dimensional subspace of $\R^3$. However there does not exist a homeomorphism of $\R^3$ which will map all of these to planes passing through the origin.\ee

We state the following corollary for the sake of completeness. 

\bc{cor2s1c5}Let $\L\subseteq \{+, 0, -\}^E$ be the set of covectors of a rank $l$ oriented matroid. Then there exists an arrangement of pseudo-hyperplanes $c\A$ such that 
\[\F(c\A) \cong (\L, \leq).\]\ec
For every $F\in \F(c\A)\setminus \{\textbf{0}\}$ let $\tilde{F}$ be the face of $\A$ such that $\tilde{F} = F\cap S^{l-1}$.  If $\sigma$ is the function assigning signs to every face then 
\begin{align*}
	\sigma(F) &:= \sigma(\tilde{F})\quad\quad \forall F\in \F(c\A)\setminus \{\mathbf{0}\}, \\
	\sigma(\mathbf{0}) &:= (0,\dots, 0).
\end{align*}
As stated before, such a face $F$ is just the cone over $\tilde{F}$, hence homeomorphic to an open polyhedral cone in $\R^l$. 

\br{rem2s1c5} In some of the literature related to topological representation theorem the term {pseudo-hyperplane} is also used for the codimension-$1$ projective space obtained by identifying antipodal points. 
Miller has used this term for topologically deformed hyperplanes in \cite{miller87} where he describes a slightly different topological representation for a certain class of oriented matroids. \er

\br{rem3s1c5}
the geometric interpretation of the composition and separation in case of pseudo-hyperplane arrangement is the same as stated in Remark \ref{rem1s1c5}.
\er 

\subsection{The microbundle complement} \label{sec2p2}
From now on we fix an arbitrary simple oriented matroid $\L$ of rank $l$, let $\A$ and $c\A$ denote the corresponding arrangements of pseudospheres (in $S^{l-1}$) and pseudo-hyperplanes (in $\R^l$) respectively. Our aim is to construct a connected subspace of $\R^{2l}$ and show that it has the homotopy type of a simplicial complex that is determined by the oriented matroid. \par

Let $c\A = \{H_1,\dots, H_n\}$ be an arrangement of pseudo-hyperplanes in $\R^l$. For every $x\in \R^l$ the \emph{arrangement localized} at $x$ is 
\[c\A_x := \{H\in c\A \mid x\in H \}. \]
Define the \emph{local complement} at $x$ as:
\[M(c\A_x) := \R^l \setminus c\A_x. \]
We now propose an analogue of the complexified complement for pseudo-hyperplane arrangements. 
Recall that for a real hyperplane arrangement $\A$ the space $M(\A)$ is the complement of the union of complexified hyperplanes inside $\C^l$. 
If one were to forget the complex structure on $\C^l$ then, topologically, it is just the tangent bundle of $\R^l$. 
The same is true for a hyperplane $H$ and its complexification $H_{\C}$. 
Hence the complexified complement of a hyperplane arrangement can also be considered as a complement inside the tangent bundle.
 
On the other hand for pseudo-hyperplane arrangements it might not be possible to define tangent spaces at every point. 
We get around this difficulty using \textit{tangent microbundles} introduced by Milnor in \cite{milnor64}.
For the benefit of the reader we begin by recalling the relevant definitions.
\bd{micro1}
A $k$-dimensional \emph{microbundle} on a topological manifold $X$ is the diagram
\[ X\stackrel{i}{\to} E\stackrel{p}{\to} X \]
consisting of the \emph{total space} $E$ and two maps, the \emph{zero section} $i: X\to E$ and the \emph{projection} $p: E\to X$. 
These are required to satisfy two properties: (1) $i$ must behave like a \emph{section}, i.e., $p\circ i = id$; and (2) $p$ must be \emph{locally trivial}, i.e., for every $x\in X$, there is a neighbourhood $V_x$ of $i(x)$ in $E$ such that the restriction $p|_{V_x}:V_x\to X$ looks like a projection $U\times \R^k \to U$.
\ed 
The difference between vector bundles and microbundles is that the main focus is near the zero-section and all requirements of linearity have been dropped.
The reader can verify that a vector bundle is indeed a microbundle; the total space and the projection maps are the same for both the bundles and the map $i$ is the zero cross-section.
It is often called the \emph{underlying microbundle} of a vector bundle.
\bd{micro2}
The \emph{tangent microbundle} of a topological manifold $X$ is defined to be the diagram
\[X\stackrel{\Delta}{\to} X\times X\stackrel{p_1}{\to} X, \]
where $\Delta$ is the diagonal map and $p_1$ is the projection on the first factor.
\ed
Intuitively it means that sufficiently close to the diagonal the fibres of $p_1$ are copies of neighborhoods of points in $X$; stacked next to each other according to their position in $X$.
If $X$ happens to be a smooth manifold then the tangent microbundle and the underlying microbundle of $TX$ are isomorphic as microbundles (see \cite[Theorem 2.2]{milnor64}).
In case of $\R^l$ since it is contractible the total spaces of two bundles are homeomorphic.

\bd{microcomplement}
The \emph{(tangent) microbundle complement} of a pseudo-hyperplane arrangement $c\A$ is the complement of the union of total spaces of the tangent microbundle of each pseudo-hyperplane in $\R^{2l}$.
We denote this space by $M(c\A)$; it can also be described as
\[M(c\A) =  \{(x, v)\in \R^l\times \R^l\mid x\in \R^l, v\in M(c\A_x)\}.\]
\ed 
This space is a natural generalization of the complexified complement since it is the complement of the family $\{H\times H\mid H\in c\A \}$ in $\R^l\times\R^l = \R^{2l}$. 

\br{remcomplement}
If the pseudo-hyperplanes of an arrangement are embedded smoothly then their tangent bundles are embedded in $\R^{2l}$.
In such a case there is one more possible generalization of the complexified complement: define it as the complement of the union of tangent bundles of pseudo-hyperplanes in $\R^{2l}$.
We call this space the \textit{tangent bundle complement}.
It will be an interesting question to figure out whether these two complements have the same homotopy type.
\er

\br{vinberg}
We should also mention that a similar construction has appeared in the literature.
Let $W$ be an infinite Coxeter group acting as a reflection group on a finite-dimensional, real vector space $V$.
Then by a theorem of Vinberg \cite{vinberg71} the action is faithful on a convex cone with non-empty interior, called the \emph{Tits' cone} and denoted $I$.
The fixed-point set of $W$ is a locally-finite collection $\A_W$ of hyperplanes in $I$.
A topological space associated with data is $I\times I\setminus \left(\bigcup_{H\in\A_W}(H\times H) \right)$. 
This space inherits the free and proper action of $W$ and the fundamental group of the quotient is an infinite-type Artin group.
\er

We now construct an open covering of the space $M(c\A)$. A connected codimension-$1$ submanifold $N$ of a manifold $M$ is \emph{two-sided} if there is a connected open neighborhood $U$ of $N$ in $M$ such that $U\setminus N$ has exactly two components each of which is open in $M$. Further, $N$ is said to be \emph{bi-collared} in $M$ if it has an open neighborhood homeomorphic to $N\times (-1, 1)$ with $N$ itself corresponding to $N\times \{0\}$. We will use the following theorem originally due to M. Brown in 1964 (see Rushing \cite[Theorem 1.7.5]{rushing73}). 

\bt{thmref3} Let $N$ be a locally flat, connected, two-sided, codimension-$1$ submanifold of $M$. Then $N$ is bi-collared in $M$. \et

Recall that \cite[Section 5]{milnor64} $N$ is said to have a \textit{microbundle neighborhood} in $M$ if there exists an open neighborhood $U$ of $N$ and a retraction $r: U\to N$ such that the diagram $N\hookrightarrow U\stackrel{r}{\to} N$ constitutes a microbundle.
This is also known as the normal microbundle.
It is clear that if $N$ has a normal microbundle inside $M$ then the embedding is locally flat.
However, the converse is true only if the submanifold is contractible.
Hence in view of \cite[Lemma 5.1]{milnor64} every face $F$ of a pseudo-hyperplane arrangement has a trivial normal microbundle.
We denote this bundle by $V_F$; note that $V_F \cong F \times (-1, 1)^{l-\dim F}$.
Choose $V_F$ such that it intersects only those faces whose closures contain $F$. 
Following lemma is now clear. 

\bl{lem3s2c5}
With the notation as above, the following statements are true:
\begin{enumerate}
\item For every $F\in \F(c\A)$ the open set $V_F$ contains $F$ and is homeomorphic to $\R^l$. 
\item If $F\leq F'$ in $\F(\A)$ then $V_F\cap V_{F'}\neq \emptyset$ and $F\nsubseteq V_{F'}$.
\item If $F$ and $F'$ are not comparable in $\F(c\A)$ then $V_F\cap V_{F'} = \emptyset$.
\end{enumerate} \el 

Consider a face-chamber pair $(F, C)$ such that closure of $C$ contains $F$. 
We denote by $C_F$ the connected component of $M(c\A_F)$ containing $C$.
Consider $\R^{2l}$ as the total space of the tangent microbundle of $\R^l$ and denote by $\tau_F$ the restriction of this bundle to $F$.
Then the total space of $\tau_F$ is $F\times \R^l$.
Use the retraction $r\colon V_F\to F$ to pullback the bundle $\tau_F$.
Let the corresponding microbundle map be $\tilde{r}\colon TV_F\to r^*(\tau_F)$ and define $W(F, C)$ to be $\tilde{r}^{-1}(F\times C_F)$.
It follows from \cite[Lemma 5.4]{milnor64} that
\[ W(F, C) = V_F \times C_F.\]

\bl{lemma1}
If $W(F, C) = W(G, D)$ then $(F, C) = (G, D)$.
\el 

\bpr
Since $W(F, C)$ is a microbundle over $V_F$, it determines $V_F$. 
By construction, $F$ is the unique minimal face contained in $V_F$; so $F = G$.

We have that $C_F = D_G$; for the sake of contradiction assume that $C\neq D$. 
But then in some neighborhood of a point $x\in F$, $C$ and $D$ must be separated by some pseudo-hyperplane $H\in c\A_x$.
So $C$ and $D$ are contained in the different connected components of $\R^l\setminus H$ which means that $C_F$ cannot the be equal to $D_G$.
\epr

\bl{lem4s2c5} $W(F_1, C_1)\cap W(F_2, C_2) \neq \emptyset$ if and only if (up to relabeling) $F_1\leq F_2$ and $C_1, C_2$ are contained in the same chamber of $\A_{F_2}$. 
\el

\bpr By construction of these open sets we have,  
\[W(F_1, C_1)\cap W(F_2, C_2) =  (V_{F_1}\cap V_{F_2})\times ((C_1)_{F_1}\cap (C_2)_{F_2}).\] 
Clearly $V_{F_1}\cap V_{F_2}\neq\emptyset$ if and only if $F_1\leq F_2$ (or $F_2\leq F_1$). 
The other intersection is non-empty if and only if, $C_1$ lies in the chamber $(C_2)_{F_2}$ of $\A_{F_2}$, i.e., $F_2\circ C_1 = C_2$ as covectors.
\end{proof}

\bd{salposet} The \emph{Salvetti poset} of a pseudo-hyperplane arrangement $c\A$ is the poset whose underlying set is $\{(F, C)\in \F \times \mathcal{C} \}$ and the ordering relation is given by 
\[(F_2, C_2) \preceq (F_1, C_1) \iff F_1\leq F_2 \hbox{~and~} F_2\circ C_1 = C_2.  \]
We denote this poset by $\mathcal{S}(c\A)$.
\ed

\br{RemSalPoset}
In the above definition we may replace the faces by the covectors and chambers by the topes of an oriented matroid.
Hence the notion of the Salvetti poset extends to non-realizable oriented matroids; such an extension first appeared in Ziegler's doctoral thesis \cite{zieg87}.
\er

Following lemma is now clear.
\bl{lemasal}
Let $(F_0, C_0), \dots, (F_k, C_k)\in \mathcal{S}(c\A)$. Then $W(F_0, C_0)\cap \dots\cap W(F_k, C_k)\neq\emptyset$ if and only if, after possibly relabelling of indices we have $(F_0, C_0)\preceq \dots\preceq (F_k, C_k)$.
\el

\bt{thm1s2c5}The collection $\{W(F, C)\mid (F, C)\in \F\times \mathcal{C}, F\leq C \}$ forms an open covering of $M(c\A)$ and whenever these open sets intersect the intersection is contractible. \et

\bpr Let $(x, v)\in M(c\A)$ be any point. 
Then there is some face $F$ such that $x\in F\subseteq V_F$ and some chamber $C$ such that $v\in C\subseteq M(c\A_x)$; hence $(x, v)\in W(F, C)$.
These sets are open and contractible because they are products of open and contractible subsets. 
Their intersections are also contractible for the same reason. \epr

Since the hypothesis of the nerve lemma \cite[Theorem 15.1]{dk1} is satisfied, the nerve of this open covering has the homotopy type of the microbundle complement $M(c\A)$. 

\bd{def1s2c5} Let $\L$ be the set of covectors of an oriented matroid the \textbf{Salvetti complex} $Sal(\L)$ is the geometric realization of the Salvetti poset. \ed

Following is the restatement of Theorem \ref{thm1s2c5} and a generalization of the classical result of Salvetti \cite{sal1}.

\bt{thm2s2c5} Let $\L$ denote the set of covectors of an oriented matroid and $c\A$ be the associated arrangement of pseudo-hyperplanes. If $M(c\A)$ is the associated space then 
\[ M(c\A) \simeq Sal(\L).\]\et

\br{remcomplement2}
In view of Remark \ref{remcomplement} one would also like to know if there exists a cell complex, described in terms of the incidence relations among the faces, that has the same homotopy type as that of the tangent bundle complement.

\er

\section{Metrical Hemisphere Complexes}\label{section3}

Recall that a central arrangement of hyperplanes decomposes the ambient Euclidean space into open polyhedral cones. As a matter of fact every hyperplane arrangement is a normal fan of a very special polytope known as the zonotope. Zonotopes can be defined in various ways: for example, projections of cubes, Minkowski sums of line segments, (polar) dual of hyperplane arrangements etc. For more on the relationship between zonotopes and hyperplane arrangements see \cite[Lecture 7]{zig95}.

\bd{defc1} A \emph{zonotope} is a polytope all of whose faces are centrally symmetric (equivalently every $2$-face is centrally symmetric). A \emph{zonotopal cell} is a (closed) $k$-cell such that its face poset is isomorphic to the face poset of a $k$-zonotope for some $k$. \ed

The face poset of a zonotope has some special combinatorial properties, the most important of which is the product structure. This product is basically the one on the face poset of a hyperplane arrangement or on the set of covectors of an oriented matroid. The following result clarifies the relationship between these three structures (see \cite[Corollary 7.17]{zig95}). 

\bt{thm1sec3} There is a natural bijection between the following three families:
\begin{enumerate}
\item the faces of a central and essential hyperplane arrangement in $\R^l$,
\item the non-empty faces of the zonotope in $\R^l$ which arises as a dual of a hyperplane arrangement,
\item signed covectors of a simple realizable oriented matroid of rank $l$.
\end{enumerate}
\et

We now extend the above correspondence to non-realizable oriented matroids. We have already shown that these oriented matroids correspond to pseudo-arrangements. In order to extend this correspondence we use the language of \emph{metrical-hemisphere complexes}. These cell complexes possess all the essential combinatorial properties of a zonotope. The metrical-hemisphere complexes (MH-complexes for short) were first introduced by Salvetti in \cite{salvetti93} to achieve a generalization of his construction and also of Deligne's theorem.\par

Let $Q$ be a connected, regular, CW complex (and $|Q|$ be the underlying space). The $1$-skeleton of such a complex $Q$ is a graph $G(Q)$ with no loops (abbreviated as $G$ if the context is clear). The vertex set of this graph will be denoted by $VG$ and the edge set by $EG$. An edge-path in $G(Q)$ is a sequence $\alpha = (l_1, \dots, l_n)$ of edges that correspond to a connected path in $|Q|$. The inverse of a path is again a path $\alpha^{-1} = (l_n,\dots, l_1)$. Two paths are composed by concatenation if the ending vertex of one of the paths is the starting vertex of another. The distance $d(v, v')$ between two vertices will be the least of the lengths of the paths joining $v$ to $v'$. Given an $i$-cell $e^i\in Q$, $Q(e^i) := \{e^j\in Q : |e^j| \subset |e^i| \}$ and $V(e^i) := VG\cap Q(e^i)$.

\bd{def2s3c3}A regular CW complex $Q$ is a QMH-complex (quasi-metrical-hemisphere complex) if there exist two maps 
\[\mino, \maxo \colon VG\times Q\to VG \] 
such that for all $v\in VG, e^i\in Q$ following properties are satisfied.
\begin{enumerate}
	\item $\mino(v, e^i)\in V(e^i)$ and $d(v, \mino(v, e^i)) = \hbox{minimum} \{d(v, u) | u \in V(e^i) \}$.
	\item $\maxo(v, e^i)\in V(e^i)$ and $d(v, \maxo(v, e^i)) = \hbox{maximum} \{d(v, u) | u \in V(e^i) \}$.
	\item $d(v, \maxo(v, e^i)) = d(v, u) + d(u, \maxo(v, e^i))$ for all $u\in V(e^i)$.
\end{enumerate}  \ed

The maps $\mino$ and $\maxo$ give us those vertices of $e^i$ which are respectively the closest and the farthest from $v$. This definition imposes a strong restriction on the $1$-skeletons of such complexes (see \cite[Proposition 1]{salvetti93}). 

\bl{lem1s3c3} If $Q$ is a QMH-complex then each circuit in $G$ has an even number of edges. Moreover if $Q$ is homeomorphic to a $k$-ball then it is a zonotopal cell. \el

For any $e^i\in Q$, indicate by $G(e^i)\subset G(Q)$ the subgraph corresponding to the $1$-skeleton of $e^i$ and by $d_{G(e^i)}$ the distance computed using $G(e^i)$. 

\bd{def3s3c3}A regular CW complex will be called an LMH-complex (local-metrical-hemisphere complex) if each $Q(e^i)$ is a QMH-complex with respect to $d_{G(e^i)}$. Moreover, the following compatibility condition also holds: 
if $e^k\in Q(e^i)\cap Q(e^j), v\in V(e^i)\cap V(e^j)$ then 
\begin{equation*} 
\mino_{(e^j)}(v, e^k) = \mino_{(e^i)}(v, e^k), \quad \maxo_{(e^j)}(v, e^k) = \maxo_{(e^i)}(v, e^k).
\end{equation*} 
Here, $\mino_{(e^j)}, \maxo_{(e^j)}$ are defined similar to $\mino, \maxo$ but using $d_{G(e^j)}$. Finally, $Q$ is an \textbf{MH-complex} if $Q$ is both a QMH-complex and an LMH-complex and for all $e^i\in Q, e^j\in Q(e^i), v\in V(e^i)$,
\begin{equation*}
\mino(v, e^j) = \mino_{(e^i)}(v, e^j), \quad \maxo(v, e^j) = \maxo_{(e^i)}(v, e^j). \end{equation*} \ed

\br{rem3s3c3}
Note that the $1$-skeleton of an MH-complex has very special properties with respect to the distance. It is not enough to have a cell complex all of whose cells are zonotopal. Here are two examples that illustrate the special nature of MH-complexes.\er
\be{mhcex1}
Consider a cell complex made up by gluing two (closed) $2$-cells, one which is octagonal and the other is trapezoidal. Three of the $1$-cells in the boundary of the trapezoidal $2$-cell are identified with three of the $1$-cells in the boundary of the octagonal cell as shown in Figure \ref{nolmh}. The resulting complex is QMH but not LMH. Consider the $1$-cell labelled by $e$ in the figure. There are two vertices, namely $v_1, v_2$, in its boundary. Considering $e$ as a member of the trapezoidal cell we see that the vertex $v_1$ is closest to the vertex $v_4$. On the other hand as a member of the octagonal $2$-cell the vertex $v_2$ is closest to $v_4$.
\begin{figure}[!ht] \begin{center}
	\includegraphics[scale=0.55,clip]{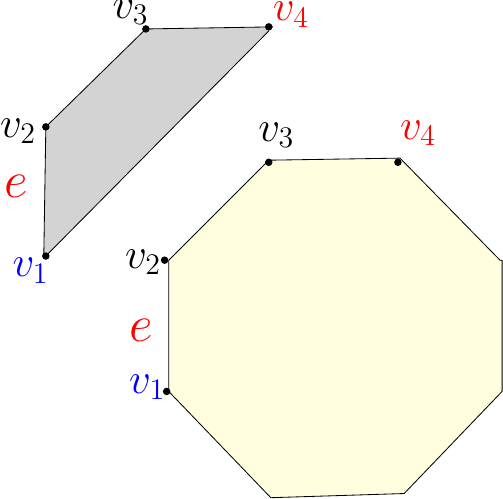}
\end{center} \caption{A QMH complex without LMH structure.} \label{nolmh}\end{figure}
\ee
\be{mhcex2}
The next example is of a cell complex obtained by replacing the trapezoidal $2$-cell from the previous complex with a $1$-cell attached as shown in Figure \ref{nomh}. The resulting cell complex is both QMH and LMH but not an MH-complex. Consider the $1$-cell labelled by $e$, it has two boundary vertices $v_1, v_2$. Considering $e$ as a member of the octagonal cell the vertex $v_2$ is closest to $v_3$. But as a member of the whole complex the boundary vertex of $e$ which is closest to $v_3$ is $v_1$.
\begin{figure}[!ht]\begin{center}
	\includegraphics[scale=0.5,clip]{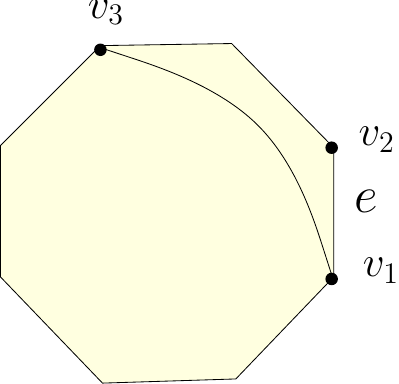}
\end{center}\caption{A QMH and LMH-complex which is not an MH-complex}\label{nomh}\end{figure}   
\ee
The following lemma establishes the combinatorial connection between zonotopes and MH-complexes. It states that the distance between any two vertices is the same no matter how it is measured, locally or globally (see \cite[Proposition 5]{salvetti93}).

\bl{lem2s3c3} Let $Q$ be an MH-complex, $e^i\in Q, v,v'\in V(e^i)$. Then $d(v, v') = d_{G(e^i)}(v, v')$. \el

\bpr Let $\alpha = (l_1,\dots, l_n)$ be a minimal path of $G(e^i)$ between $v$ and $v'$ (so $d_{G(e^i)}(v, v') = n$). Let $v_{j-1}, v_j$ be the vertices of $l_j$ ordered according to the orientation of $\alpha$ from $v$ to $v'$. Since $\alpha$ is minimal in  $G(e^i)$ and $Q$ is an MH-complex we have,
\[\mino_{(e^i)}(v, l_j) = v_{j-1} = \mino(v, l_j). \]
Hence, $d_{G(e^i)}(v, v_j) = d_{G(e^i)}(v, v_{j-1}) + 1$ and $d(v, v_j) = d(v, v_{j-1}) + 1$ for $j = 1,\dots, n$.
The statement then follows by induction. \epr

Given a pseudo-hyperplane arrangement  $c\A$ in $\R^l$ (corresponding to an oriented matroid) and $\F(c\A)$ as its face poset let $(\R^l, \F^*)$  we denote the cell complex, embedded in $\R^l$, which is dual to the induced stratification. 
Where $\F^*$ is the face poset of the dual stratification. 
Note that the $\F^*$ is the dual face poset, i.e., the poset obtained by reversing the order in $\F$. 
For every chamber of $c\A$ there corresponds a $0$-cell in $(\R^l, \F^*)$ and in general for a $k$-face $F$ there is an $(l-k)$-cell $F^*$. 
We now state the theorem that generalizes the relationship between hyperplane arrangements and zonotopes. 

\bl{lem1s2c5} Let $c\A$ be a pseudo-hyperplane arrangement corresponding to an oriented matroid. Then the dual cell complex $(\R^l, \F^*)$ is an MH-complex. \el

\bpr 
First we need to define the two maps $\mino$ and $\maxo$.
For notational simplicity we will not differentiate between a face and its dual.
Given a face $F$ and a chamber $C$ define $\mino(F, C)$ to be the chamber incident with $F$ which is closest to $C$. 
Note that in view of Remark \ref{rem3s1c5} such a chamber always exists and is unique.
Analogously one can define $\maxo(F, C)$ to be the unique chamber incident with $F$ which is farthest from $C$.
Recall that a path between two chambers (or vertices) $C, C'$ is minimal if and only if it crosses each separating pseudo-hyperplanes exactly once and does not cross any other any other.
By abuse of notation we do not differentiate between the separation set of two topes and the set of pseudo-hyperplanes separating the corresponding chambers which we denote by $S(C, C')$.
Note that given a face $F$ such that $F = H_1\cap\cdots\cap H_k$ then $S(\mino(F, C), \maxo(F, C)) = \{H_1,\dots, H_k\}$. In addition, if $D$ is a chamber incident with $F$ then 
\[S(\mino(F, C), D)\cup S(\maxo(F, C), D) =  \{H_1,\dots, H_k\}\hbox{~and~} S(\mino(F, C), D)\cap S(\maxo(F, C), D) =\emptyset.\] 
The distance between two chambers $C, C'$ is defined as the cardinality of the corresponding separating set.
Now the following equalities establish that $(\R^l, \F^*)$ is a QMH-complex.
\[ S(C, D) = S(C, \mino(F, C))\cup S(\mino(F, C), D),\quad S(C, \mino(F, C))\cap S(\mino(F, C), D)=\emptyset.\]
\[ S(C, \maxo(F, C)) = S(C, D)\cup S(D, \maxo(F, C)), \quad S(C, D)\cap S(D, \maxo(F, C))=\emptyset.\]
The proof that $(\R^l, \F^*)$ is an MH-complex follows from a straightforward verification of the remaining axioms.\epr

We should note here that the MH-complex structures on the closed unit ball in $\R^l$ are in one-to-one correspondence with simple oriented matroids of rank $l$. This can be done by proving that the cells of this MH-complex satisfy covector axioms. 

For the sake of completeness we will explicitly describe the cells of the Salvetti complex. The $0$-cells correspond to topes of the oriented matroid which we denote by $[T, T]$. Let $X$  be a covector which corresponds to an $(l-k)$-face $F_X$ of the corresponding pseudo-hyperplane arrangement. For every such covector $X$ and every tope $T$ such that $X\leq T$ there corresponds a $k$-cell $[X, T]$ which is homeomorphic to $F_X^*$. The boundary of such a cell is given by:
\[\partial [X, T] = \bigcup_{X < Y}[Y, Y\circ T]. \]
The Salvetti complex has an oriented $1$-skeleton by directing an edge $[X, T]$ from $[T, T]$ to $[T', T']$, where $T$ and $T'$ are topes greater than $X$.  \par
The oriented $1$-skeleton of a Salvetti complex has an additional combinatorial structure which we now describe. Recall that oriented matroids can also be defined using the \emph{circuit} axioms (see \cite[Definition 3.2.1]{ombook99} and \cite[Section II]{fl78}).
\bd{defcc}
An oriented matroid is a triple $(E, \mathcal{C}, \ast)$, where $E$ is a finite set, $\mathcal{C}$ is a collection of non-empty subsets of $E$, and $\ast$ is fixed-point free involution on $E$, such that:
\begin{enumerate}
	\item if $A, B\in \mathcal{C}$ with $A\subset B$, then $A = B$;
	\item if $A\in\mathcal{C}$ then $A^{\ast}\in \mathcal{C}$ and $A\cap A^{\ast} = \emptyset$;
	\item for all $A, B\in \mathcal{C}$, $x\in A\cap B^{\ast}$, and $A\neq B^{\ast}$ there is a $D\in\mathcal{C}$ with $D\subseteq (A\cup B)\setminus \{x, x^{\ast}\}$.
\end{enumerate}
\ed
The following result is a straightforward verification of the above definition and hence left to the reader as an exercise. 

\bp{lemlastsec3} Let $Sal(\L)_1$ denote the oriented $1$-skeleton of a Salvetti complex associated to an oriented matroid $\L$. Define a map $\phi$ on the closed $1$-cells by sending $[X, T]\mapsto [X, T']$ where $T, T'$ are topes covering $X$. If $\mathcal{C}$ is the set of all circuits of the graph underlying $Sal(\L)_1$ then the triple $(Sal(\L)_1, \mathcal{C}, \phi)$ satisfies the circuit axioms of an oriented matroid. \ep

 \section{Topology of the microbundle complement}\label{section4}

In this section we further investigate the interaction between the combinatorics of the oriented matroid and the topology of the space $M(c\A)$. The results in this section are not necessarily new. Most of them are already known in the context of hyperplane arrangements. We state them in the general context of oriented matroids; also note that the proofs adapt in this situation. 

\subsection{Salvetti complex} We start by elaborating on the relationship between an oriented matroid and the combinatorics of the cells of the associated Salvetti complex. 

\bt{thm1sec4}Let $\L$ be an oriented matroid and $Sal(\L)$ be the associated Salvetti complex. Then we have the following:
\begin{enumerate}
	\item The complex $Sal(\L)$ is an MH-complex.
	\item The number of $0$-cells of $Sal(\L)$ is equal to the number of its $l$-cells which is also equal to the number of topes of $\L$.
	\item Every chain in the barycentric subdivision of $Sal(\L)$ corresponds to pair containing a chain in $(\L, \leq)$ and a tope.
	\item The geometric realization of $(\L, \leq)$ is a retract of $Sal(\L)$. 
	\item $\chi(Sal(\L)) = 0$.
	\item The homeomorphism type of $M(c\A)$ is completely determined by $\L$. 
\end{enumerate} \et

\bpr 
The first statement follows from the construction of the Salvetti complex. The $1$-skeleton of this complex is obtained by doubling the edges in the $1$-skeleton of $(\R^l, \F^*)$. The second statement follows from Definition \ref{def1s2c5}, since the $0$-cells correspond to tope-tope pairs. The top-dimensional cells correspond to pairs of the form $[\mathbf{0}, T]$ for every tope $T$.\par 
The proof of the third statement is also clear from the definition. As for statement (4) the map $X\mapsto [X, X\circ T]$ is an inclusion of $\Delta(\F(\L))$ into $Sal(\L)$ and the map in the other direction defined by $[X, T]\mapsto X$ is a retraction.\par 
For (5) one has to take the difference of the number of odd-dimensional cells and the number of even-dimensional cells of $(\R^l, \F^*)$. This can be done by analyzing the link of every vertex; see \cite[Theorem 3.3.6]{deshpande_thesis11} for details. This statement also implies that there are no bounded chambers. The last statement follows from Bj\"orner-Ziegler \cite[Theorem 5.3]{bz92}. \epr

Recall that a path in a regular cell complex is a sequence of consecutive edges and its length is the number of edges. A minimal path is a path of shortest length among all the paths joining its end points. In case of an oriented $1$-skeleton, by a positive path we mean a path all of whose edges are traversed in the positive direction. 
The next couple of results are about the fundamental groupoid of the Salvetti complex. As in the case of hyperplane arrangements a directed path in the Salvetti complex is basically a sequence of adjacent chambers. Moreover, a path is minimal if and only if it does not cross a pseudo-hyperplane twice.

\bl{lem5s2c5} With the notation as before, any two positive minimal paths in the $1$-skeleton of $Sal(\L)$ that have same initial as well as terminal vertex are homotopic relative to $\{0, 1\}$. \el

\bpr Given two positive minimal paths $\alpha, \beta$ in $Sal(\L)$ with the same end points apply the retraction map to get paths in $\Delta(\F^*)$. Observe that no two edges of these two paths are sent to a same edge of $\Delta(\F^*)$. The conclusion follows from the fact that the geometric realization $\Delta(\F^*)$ is contractible.\epr

For a combinatorial proof see Cordovil \cite[Theorem 2.4]{cordovil93}.\par 

Given an oriented matroid $\L$ let $\G^+$ denote the associated \emph{positive category}. 
It is the category of directed paths in $Sal(\L)$. 
The objects of this category are the vertices of the Salvetti complex and morphisms are directed homotopy classes of positive paths (i.e., two such paths are equivalent if and only if they are connected by a sequence of substitutions of minimal positive paths with the same end points). 
Let $\G$ denote the fundamental groupoid of the associated Salvetti complex. Note that it is also the category of fractions of $\G^+$. We now show that the cancellation law holds in $\G^+$; this is a well-known fact in case of hyperplane arrangements (see for example \cite[Proposition 1.12]{deli72}).

\bt{thm2sec4} The associated canonical functor $J\colon \G^+\to \G$ is faithful on the class of minimal positive paths.\et

\bpr We have seen that any two minimal positive paths with the same end points represent the same morphism in $\G$. We now have to show that they also represent the same morphism in $\G^+$. Given two such paths we will show that one of them can be obtained from the other using a sequence of substitutions of minimal positive paths. \par
Let $\alpha, \beta$ be two minimal positive paths in $Sal(\L)$ from a vertex $[C, C]$ to another vertex $[D, D]$. For notational simplicity we will not differentiate between a chamber and the corresponding vertex in the Salvetti complex, i.e., we use $C$ in place of $[C, C]$. 
We represent these paths by edge sequences, i.e., $\alpha = (a_1,\dots, a_n)$ and $\beta = (b_1,\dots, b_n)$. We proceed by induction on $n$, the cases $n = 0, 1$ being trivial. Assume that the statement is true for all positive minimal paths of length $n-1$. Now, if $a_1 = b_1$ then $(a_2,\dots, a_n)$ is equal to $(b_2,\dots, b_n)$ in $\G$, by Lemma \ref{lem5s2c5} so that $\alpha$ and $\beta$ represent the same morphism in $\G^+$.\par 

If $a_1\neq b_1$ let $H_a$ and $H_b$ be the pseudo-hyperplanes dual to $a_1$ and $b_1$ respectively. For $H\in \A$ that does not separate the chambers $C, D$ let $X_H(C, D)$ denote the closure of the connected component of $\R^l\setminus H$ that contains both $C, D$. Let $\h(C, D)$ be the intersection of all such $X_H(C, D)$'s. Then it is clear that either $\h(C, D) = \emptyset$ (this happens when $C$ and $D$ are opposite chambers) or $\h(C, D)$ is connected and homeomorphic to the intersection of finitely many closed half-spaces. Let us first look at the sub-case when $\h(C, D)\neq \emptyset$. Then the intersection $\h(C, D)\cap H_a\cap H_b$ contains a codimension-$2$ face, say $F^2$. 
Then there is a positive minimal path from $C' := \maxo(C, F^2)$ which we denote by $\gamma_0$. Let $\gamma_1, \gamma_2$ be two positive minimal paths starting from $\maxo(C, a_1)$ and $\maxo(C, b_1)$ and both ending at $C'$. Then $\alpha$ is equivalent to $a_1\gamma_1\gamma_0$ because they share the same first edge, and they have the same beginning and ending points. For the same reasons $\beta$ and $b_1\gamma_2\gamma_0$ are equivalent. However, it is not hard to show that $a_1\gamma_1\gamma_0$ and $b_1\gamma_1\gamma_0$ are also equivalent, hence concluding the proof in this sub-case. The other sub-case when $\h(C, D)$ is empty can be treated similarly.\epr

Delucchi in his thesis \cite{del1} introduced the theory of Salvetti-type diagram models in 
order to characterize and classify covering spaces of the complexified complement of a hyperplane arrangement. This homotopy theoretic 
technique also works in case of non-realizable oriented matroids. The covering spaces can be realized as the homotopy colimit of certain 
diagrams of spaces defined using covering groupoids of $\G$. Stating these results requires some terminology from homotopy theory and 
would be a digression hence we refer the reader to \cite[Chapters 3, 5]{deshpande_thesis11} for precise statements. 

\subsection{Cohomology of the complement}The cohomology algebra of the complexified complement of a hyperplane arrangement is known as the Orlik-Solomon algebra (OS-algebra) and is determined by the intersection lattice. 
The construction of this OS-algebra is completely combinatorial and works for any matroid (see \cite[Chapter 3]{orlik92}). 
Hence even in the case of a non-realizable oriented matroid its underlying matroid has an associated OS-algebra. 
It was shown by Gel'fand and Rybnikov \cite[Theorem 5]{gr89} that the cellular cohomology ring of the Salvetti complex is isomorphic to the associated OS-algebra. 
Therefore the Salvetti complexes have the `right' cohomology. 
Their technique was generalized by Bj\"orner and Ziegler \cite[Corollary 7.3]{bz92} to arbitrary complex arrangements to give a completely combinatorial proof of the Brieskorn-Orlik-Solomon theorem. \par

It was proved by Dimca and Papadima in \cite{dimca_papadima03} that the complement of a hyperplane arrangement is minimal, i.e., it has the homotopy type of a CW complex such that number of $k$-cells is equal to $k$-th Betti number. 
From the work of Delucchi \cite{delucchi08} it follows that in case of a non-realizable oriented matroid the associated space $M(c\A)$ is minimal. 
This is done by constructing maximum acyclic matchings of the face poset of the Salvetti complex such that its critical cells are in one-to-one correspondence with the topes, this correspondence is achieved via no-broken-circuits \cite[Proposition 2, Lemma 5.10]{delucchi08}. 

\subsection{Simplicial oriented matroids} We now turn to simplicial arrangements, that is, arrangements in which every chamber is a cone over an open simplex. Alternately, an oriented matroid is simplicial if $\L\setminus \{\textbf{0}\}$ is isomorphic to the face poset of a simplicial decomposition of the sphere (or for every tope $T$ the interval $[\textbf{0}, T]$ is Boolean). \par
Recall that there is an ordering on the topes of an oriented matroid induced by the corresponding chambers. 
Fix a tope $T$. For any other tope $S$ define the distance $d(S, T)$ between $S$ and $T$ to be the number of pseudo-hyperplanes that separate the corresponding chambers. 
Now, for another tope $S'$,  $S'\preceq_T S$ if and only if $d(S', T) \leq d(S, T)$ is a partial order. Denote this poset by $\P_T(\L)$.
 
\bt{thm5s2c5} Let $\L$ be an oriented matroid and $c\A$ be the associated pseudo-hyperplane arrangement. Then the following are equivalent.
\begin{enumerate}
	\item $\L$ is simplicial.
	\item The positive category $\G^+$ admits the Deligne normal form.
	\item The tope poset $\P_T(\L)$ is a lattice for every tope $T$.
\end{enumerate} All of the above conditions imply that the space $M(c\A)$ is $K(\pi, 1)$. \et

\bpr The Deligne normal form is a particular factorization of loops in the Salvetti complex into compositions of positive paths. We refer the interested reader to Paris \cite{pa1} and Delucchi \cite[Chapter 6]{del1} for an insight into Deligne's original arguments.\par 
(1) $\Rightarrow$ (2) is originally due to Deligne \cite{deli72} (and a reproof by Paris \cite{pa1}); both these proofs are for realizable oriented matroids. For non-realizable oriented matroids the proof was given by Cordovil \cite[Theorem 4.1]{cordovil94b} and by Salvetti \cite[Theorem 33]{salvetti93}. 2 $\Rightarrow$ 1 is due to Paris \cite{pa4}. For 1 $\iff$ 3 see 
\cite{bez90} and the proof of 3 $\iff$ 1 is in Delucchi's thesis \cite[Theorem 6.4.6, Lemma 6.5.2]{del1}. \epr

It will be interesting to see which properties of the pure Artin groups are also satisfied by the fundamental groups of Salvetti complexes associated with simplicial oriented matroids. For example, one would like to see if Charney's arguments in \cite{charney92} generalize in this context to show that the fundamental group of the Salvetti complex is bi-automatic or in general is an example of a Garside group.\par 

Obvious examples of arrangements that were not covered by Deligne's theorem are the simplicial arrangements of pseudolines. 
A simplicial arrangement of pseudolines in $\R P^2$ consists of a finite family of simple closed curves such that every two curves have precisely one point in common and every $2$-face is isomorphic to a triangle. 
By applying the coning process we get an arrangement of (non-stretchable) pseudoplanes in $\R^3$ whose face poset corresponds to a rank $3$ non-realizable oriented matroid. 
The Salvetti complex associated to such oriented matroids is a $K(\pi, 1)$ space. 
In fact at least seven infinite families of non-stretchable simplicial arrangements of pseudolines are known, see \cite[Chapter 3]{grunbaum67} for details and examples.\par 

We end the paper by a discussion on the relevance of our work. 
We wish to state here that the topology of the Salvetti complex is different in the realizable and non-realizable case even though the above results suggest otherwise. 
The main evidence is the recent work of Delucchi and Falk \cite{delfalk13}. 
The authors show that the fundamental group of the Salvetti complex associated to the non-Pappus arrangement is not isomorphic to any realizable arrangement group. 
Hence there is a lot of scope to uncover various combinatorial aspects of non-realizable oriented matroids that govern the topology of the Salvetti complex. 
For example, it is not known whether in the non-realizable case the cohomology algebra of the Salvetti complex is formal in the sense of Sullivan. 
A negative answer to this question will provide more insight into oriented matroids whereas an affirmative answer will provide a completely combinatorial proof of formality of arrangement complements.

\renewcommand{\bibname}{References}
\bibliographystyle{abbrv} 
\bibliography{arrman} 

\end{document}